\documentclass[3p]{elsarticle}

\usepackage{amsmath}
\usepackage{amsthm}
\usepackage{amssymb}
\usepackage{bm}
\usepackage{accents}

\newtheorem{thm}{Theorem}

\newdefinition{rmk}{Remark}
\newproof{pf}{Proof}

\newcommand{\grad}{\mathop{\rm grad}\nolimits}
\newcommand{\const}{\mathop{\rm const}\nolimits}
\renewcommand{\div}{\mathop{\rm div}\nolimits}

\numberwithin{equation}{section}

\journal{arXiv} 

\graphicspath{{figs/}}

\begin{document}

\begin{frontmatter}

\title{Numerical solution of the heat conduction problem with memory\tnoteref{label1}}
\tnotetext[label1]{The publication was financially supported by the research grant 20-01-00207 of Russian Foundation for Basic Research}

\author{P.N. Vabishchevich\corref{cor1}\fnref{lab1,lab2}}
\ead{vabishchevich@gmail.com}
\cortext[cor1]{Correspondibg author.}

\address[lab1]{Nuclear Safety Institute, Russian Academy of Sciences,
              52, B. Tulskaya, 115191 Moscow, Russia}

\address[lab2]{North-Caucasus Center for Mathematical Research, North-Caucasus Federal University, 
                   1, Pushkin Street, 355017 Stavropol, Russia}

\begin{abstract}
It is necessary to use more general models than the classical Fourier heat conduction law to describe small-scale thermal conductivity processes.
The effects of heat flow memory and heat capacity memory (internal energy) in solids are considered in first-order integrodifferential evolutionary equations with difference-type kernels.
The main difficulties in applying such nonlocal in-time mathematical models are associated with the need to work with a solution throughout the entire history of the process.
The paper develops an approach to transforming a nonlocal problem into a computationally simpler local problem for a system of first-order evolution equations. Such a transition is applicable for heat conduction problems with memory if the relaxation functions of the heat flux and heat capacity are represented as a sum of exponentials. The correctness of the auxiliary linear problem is ensured by the obtained estimates of the stability of the solution concerning the initial data and the right-hand side in the corresponding Hilbert spaces.
The study's main result is to prove the unconditional stability of the proposed two-level scheme with weights for the evolutionary system of equations for modeling heat conduction in solid media with memory. In this case, finding an approximate solution on a new level in time is not more complicated than the classical heat equation.
The numerical solution of a model one-dimensional in space heat conduction problem with memory effects is presented.
\end{abstract}


\end{frontmatter}
\section{Introduction}\label{s:1}

Classical thermal conduction is based on the Fourier heat conduction law, where the heat flux is proportional to the temperature gradient.
In this approximation, the rate of heat propagation is infinite.
When considering heat transfer processes on small scales, we are forced to use more complex heat transfer models \cite{joseph1989heat,straughan2011heat}. 
A classic example of non-Fourier heat conduction equations is the hyperbolic heat conduction model when the effects of heat flux relaxation are taken into account. 
The mathematical model becomes more complicated but remains local: we have a second-order evolution equation for temperature or a system of first-order equations for temperature and heat flux.
In the general case, we should focus on nonlocal heat transfer models, where the thermal state depends on prehistory.
Memory effects of heat flow and heat capacity are accounted for using nonlocal evolution equations for temperature \cite{gurtin1968general,nunziato1971heat} with integral terms with difference-type kernels.

The problems of numerical solution of nonlocal in time IBVPs are well known. 
We can use some or other quadrature formulas to approximate the integral time terms and have suitable properties on stability and convergence of the approximate solution to the exact one (see, for example, \cite{ChenBook1998}).
The main problems are related to the computational cost of finding the solution on the new time level.
We have to operate with approximate solutions for all prior time moments. 
Because of this, we are particularly interested in approaches that allow us to go from a nonlocal model to a local model not much more complex than the classical heat conduction equation.

We use in our work a well-known method of approximate solution of the Volterra integral equations \cite{linz1985analytical}, which is associated with an approximation of the kernel by a sum of exponents. 
The construction of such approximations is an independent nonlinear approximation problem \cite{braess1986nonlinear}.
The justification of such a transition for evolutionary first-order integrodifferential equations with solution memory and solution derivative in time is given in our works \cite{vabMemory,vabMemoryB}. 
The present paper considers evolution equations for heat conduction problems with heat flow memory and heat capacity memory. 

The paper is organized as follows. In Section~2, we describe nonlocal heat conduction problems with heat flow memory and heat capacity memory. In Section~3, we discuss the transition from a nonlocal problem with heat flow and energy distribution functions as a sum of exponents to the Cauchy problem for a local system of first-order evolution equations.
Section~4 is devoted to the construction and study of two-level operator-difference schemes. 
In Section~5, we present examples of numerical solutions of one-dimensional spatial heat conduction problems with memory.
The final section is the conclusion.

\section{Nonlocal heat conduction problems}\label{s:2}

We limit ourselves to considering linear models of heat conduction of a solid isotropic body.
In classical heat conduction theory, the heat flow $\bm q(\bm x,t)$ at a point $\bm x$ at some time $t$ is determined by the temperature gradient $T(\bm x,t)$ (Fourier's law of heat conduction) when
\begin{equation}\label{2.1}
 \bm q = - k \grad T ,
\end{equation} 
where $k = \const$ is the thermal conductivity. 
For a rigid conductor the internal energy equation is
\begin{equation}\label{2.2}
 \frac{\partial e}{\partial t} + \div \bm q = r ,
\end{equation} 
where $r(\bm x,t)$ is the distributed heat suction.
The energy increments $e(\bm x,t)$ are due to changes in temperature, so that
\begin{equation}\label{2.3}
 e = c \, T + \const, 
\end{equation} 
where $c = \const$ is the heat capacity.
Substituting (\ref{2.1}), (\ref{2.3}) into (\ref{2.2}) gives the heat conduction equation
\begin{equation}\label{2.4}
 c \frac{\partial T}{\partial t} - k \div \grad T = r .
\end{equation} 

Mathematical models of heat conduction with memory are proposed in \cite{gurtin1968general,nunziato1971heat}. 
The prehistory of the process is accounted in \cite{gurtin1968general} by specifying the heat flow by the relation
\begin{equation}\label{2.5}
 \bm q(\bm x,t) = - \int_{0}^{t} \alpha (t-s) \grad T (\bm x,s) ds ,
\end{equation} 
where $\alpha (t)$ is a positive, decreasing function relaxation of the heat flux.
In the special case $\alpha (t) = k \delta (t)$ ($\delta (t)$ is a $\delta$-function), the Fourier law follows from  (\ref{2.1}), (\ref{2.5}).
The relation (\ref{2.5}) leads us instead of the equation (\ref{2.4}) to the evolutionary integrodifferential equation
\begin{equation}\label{2.6}
 c \frac{\partial T}{\partial t} - \int_{0}^{t} \alpha (t-s) \div \grad T (\bm x,s) ds = r . 
\end{equation} 

Let's assume that 
\begin{equation}\label{2.7}
 \alpha (t) = \frac{k}{\tau_q} \exp \Big (- \frac{t}{\tau_q} \Big ) ,
\end{equation}  
where $\tau_q$ is the relaxation time. In this case, from (\ref{2.5}) we have the Maxwell-Cattaneo law \cite{cattaneo1948sulla} 
\begin{equation}\label{2.8}
 \tau_q \frac{\partial \bm q}{\partial t} + \bm q = - k \grad T 
\end{equation} 
for the fluxes. This leads us to the hyperbolic heat conduction equation \cite{bubnov1976wave} 
\begin{equation}\label{2.9}
 c \tau_q \frac{\partial^2 T}{\partial t^2} + c \frac{\partial T}{\partial t} - k \div \grad T 
 = r + \tau_q \frac{\partial r}{\partial t}.
\end{equation} 
The fundamental point is that instead of a nonlocal model (\ref{2.5}) for the exponential dependence of the relaxation function on time (\ref{2.7}), we have a local model (at each individual point $(\bm x,t)$) for the dependence of heat flux on temperature gradient as a differential relation (\ref{2.8}). Instead of the integrodifferential equation (\ref{2.6}), a second-order evolution equation (\ref{2.9}) is used.

For the heat flux relaxation function we will use \cite{nunziato1971heat} the expression
\begin{equation}\label{2.10}
 \alpha(t) = k \delta (t) + \tilde{\alpha} (t) .
\end{equation}  
Similar relations are used for internal energy relaxation when
\begin{equation}\label{2.11}
 e(\bm x,t) = \int_{0}^{t} \beta (t-s) T (\bm x,s) ds + \const, 
\end{equation} 
\begin{equation}\label{2.12}
 \beta (t) = c \delta (t) + \tilde{\beta} (t) .
\end{equation} 
From (\ref{2.2}), (\ref{2.5}), (\ref{2.10})--(\ref{2.12}), we obtain the heat conduction equation with memory
\begin{equation}\label{2.13}
\begin{split}
 c \frac{\partial T}{\partial t} & + \frac{\partial }{\partial t} \int_{0}^{t} \tilde{\beta} (t-s) T (\bm x,s) ds \\
 & - k \div \grad T - \int_{0}^{t} \tilde{\alpha} (t-s) \div \grad T (\bm x,s) ds = r . 
\end{split}
\end{equation} 

We consider a numerical solution of the boundary value problem for the equation (\ref{2.13}).
For simplicity, we restrict ourselves to the case of a bounded domain $\Omega$ with a boundary $\partial \Omega$ on which homogeneous Dirichlet boundary conditions are given:
\begin{equation}\label{2.14}
 T (\bm x,t) = 0,
 \quad \bm x \in \partial \Omega .
\end{equation} 
The initial thermal state is defined according to
\begin{equation}\label{2.15}
 T (\bm x,0) = T_0 (\bm x), 
 \quad \bm x \in \partial \Omega .
\end{equation} 
Single-valued solvability of boundary value problems of type (\ref{2.13})--(\ref{2.15}) in the corresponding functional classes, the qualitative properties of the solutions are discussed, for example, in \cite{davis1976hyperbolicity,davis1978linear}.

The object of our consideration is computational algorithms for approximate solutions of heat conduction problems with allowance for memory effects. The main problems are generated by the nonlocality of the equation (\ref{2.13}), the fact that the solution depends on the complete prehistory of the process.
When finding the solution at a new level in time, we have to consider the solution at previous levels in time, which, in particular, is associated with high computational and memory costs.
Distinct possibilities of a principal reduction of computational work are seen on the example of the transition to the local equation (\ref{2.9}) due to transition from the integrodifferential relation (\ref{2.5}) to the differential relation (\ref{2.8}) when the heat flow distribution function is given as (\ref{2.7}). 
We can expect a similar simplification of the problem when the heat flux and energy distribution functions are represented as
\begin{equation}\label{2.16}
 \tilde{\alpha}(t) = \sum_{i = 1}^{m} \alpha_i \exp(- \mu_i t),
 \quad \tilde{\beta }(t) = \sum_{j = 1}^{l} \beta_j \exp(- \nu_j t) . 
\end{equation} 
For general dependencies $\tilde{\alpha}(t), \tilde{\beta }(t)$ this representation corresponds to using nonlinear approximations of functions by the sum of exponents \cite{braess1986nonlinear}.  

\section{Local system of equations}\label{s:3}

After a finite-element or finite-volume approximation over the space, we consider the discrete problem in the corresponding finite-dimensional Hilbert space $H$ \cite{KnabnerAngermann2003,QuarteroniValli1994}.
For the scalar product and norm we use the notation
$(u, v)$ and $\|u\|=(u,u)^{1/2}$ for $u, v \in H$ respectively.
A positively determined selfadjoint operator $D$ is associated with the Hilbert space $H_D$, in which $(u, v)_D = (Du, v)$, $\|u\|_D = (D u, u)^{1/2}$. 

The dimensionless problem (\ref{2.13})--(\ref{2.15}) leads us to the equation
\begin{equation}\label{3.1}
 \frac{d u}{d t} + \frac{d }{d t} \int_{0}^{t} \tilde{\beta} (t-s) u(s) ds 
 + A u + \int_{0}^{t} \tilde{\alpha} (t-s) A u(s) ds = f(t) ,
 \quad t > 0 , 
\end{equation} 
for $u(\cdot, t) = u(t)$.
The operator $A$ is associated with an approximation of the Laplace operator for functions satisfying boundary conditions (\ref{2.14}).
It is natural to assume that in $H$ it is constant (independent of $t$) and $A = A^* > 0$. 
We consider the Cauchy problem for the Volterra evolutionary equation (\ref{3.1}) when we use the initial condition
\begin{equation}\label{3.2}
 u(0) = u_0 .
\end{equation} 

The transition from a nonlocal to a local problem when approximating the difference-type kernel by a sum of exponents is provided by introducing new desired variables.
For the problem (\ref{2.16}), (\ref{3.1}), (\ref{3.2}), we put
\[
 v_i(t) = \int_{0}^{t} \exp(- \mu_i (t-s)) u(s) ds ,
 \quad i = 1,2, \ldots, m, 
\]
\[
 w_j(t) = \int_{0}^{t} \exp(- \nu_j (t-s)) u(s) ds ,
 \quad j = 1,2, \ldots, l. 
\]
For these auxiliary functions, we have the local relations
\begin{equation}\label{3.3}
 \frac{d v_i}{d t} + \mu_i v_i - u(t) = 0 ,
 \quad i = 1,2, \ldots, m, 
\end{equation} 
\begin{equation}\label{3.4}
 \frac{d w_j}{d t} + \nu_j w_j - u(t) = 0 ,
 \quad j = 1,2, \ldots, l. 
\end{equation} 
Equation (\ref{3.1}) is written as
\[
 \frac{d u}{d t} + \sum_{j=1}^{l} \beta_j \frac{d w_j}{d t} 
 + A u + \sum_{i=1}^{m} \alpha_i A v_i = f(t) .
\]
Given (\ref{3.4}), we have
\begin{equation}\label{3.5}
 \frac{d u}{d t} + \sum_{j=1}^{l} \beta_j u + A u - \sum_{j=1}^{l} \beta_j \nu_j w_j
 + \sum_{i=1}^{m} \alpha_i A v_i = f(t) .
\end{equation}  
The initial condition (\ref{3.2}) for $u(t)$ is supplemented by initial conditions for auxiliary functions:
\begin{equation}\label{3.6}
 v_i(0) = 0 ,
 \quad i = 1,2, \ldots, m, 
 \quad w_j(0) = 0,
 \quad j = 1,2, \ldots, l.  
\end{equation} 
Thus from the nonlocal problem (\ref{2.16}), (\ref{3.1}), (\ref{3.2}) for $u(t)$, we arrive at the local Cauchy problem for the system of equations (\ref{3.2})--(\ref{3.6}).

The main result concerns the stability for the initial data and the right-hand side for the system of equations (\ref{3.3})--(\ref{3.5}).
The correctness of this linear problem is established under the assumption that the coefficients in the representations (\ref{2.16}) are positive:
\begin{equation}\label{3.7}
 \alpha_i > 0 ,
 \quad \mu_i > 0 ,
 \quad i = 1,2, \ldots, m, 
 \quad \beta_j > 0,
 \quad \nu_j > 0,
 \quad j = 1,2, \ldots, l.   
\end{equation}  

\begin{thm}\label{t-1}
If (\ref{3.7}) is satisfied, for the solution of the Cauchy problem (\ref{3.2})--(\ref{3.6}), there is an a priori estimate 
\begin{equation}\label{3.8}
 \|y(t)\|_* \leq \|u_0\| + \int_{0}^{t} \|f(s\| ds ,
\end{equation} 
where
\[
 \|y(t)\|_*^2 = \|u(t)\|^2 + \sum_{i=1}^{m} \alpha_i \|v_i(t)\|_A^2 + \sum_{j=1}^{l} \beta_j \nu_j \|w_j(t)\|^2 .
\] 
\end{thm}

\begin{pf} 
Let us multiply equation (\ref{3.5}) scalarly in $H$ by $u(t)$, separate equations (\ref{3.3}) by $\alpha_i A v_i(t)$, and equations (\ref{3.4}) by $\beta_j \nu_j w_j(t)$. The summation of these equations gives
\begin{equation}\label{3.9}
\begin{split}
 \frac{1}{2} & \frac{d}{d t} \Big (\|u(t)\|^2 + \sum_{i=1}^{m} \alpha_i \|v_i(t)\|_A^2 + \sum_{j=1}^{l} \beta_j \nu_j \|w_j|^2 \Big ) 
 + \|u(t)\|_A^2 \\
 & + \sum_{i=1}^{m} \alpha_i \mu_i \||v_i(t)\|_A^2 + \sum_{j=1}^{l} \beta_j \|u(t) - \nu_j w_j\|^2 = (f,u) .
\end{split}
\end{equation} 
Given 
\[
 (f,u) \leq \|f(t)\| \|u\\| \leq \|f(t)\| \||y(t)\|_* ,
\] 
from (\ref{3.9}), we have the desired estimate (\ref{3.8}). 
\end{pf}

From (\ref{3.8}), a simpler estimate follows:
\[
 \|u(t)\| \leq \|u_0\| + \int_{0}^{t} \|f(s\| ds . 
\]
This a priori estimate is standard for the problem without memory effects ($\tilde{\alpha} (t) = \tilde{\beta} (t) \equiv 0$ in equation (\ref{3.1})).

\section{Two-level difference scheme}\label{s:4}

When approximating first-order evolution equations, it is natural to focus on standard two-level schemes with weights \cite{LeVeque2007,SamarskiiTheory}. Without limiting generality, we will assume that the time grid is uniform, i.e.
$t^n = n \tau, \ n = 0, 1, \ldots,$ where $\tau > 0$ is time step. 
We denote the approximate solution at time $t^n$ by $u^n, \ \ n = 0, 1, \ldots$ and let
\[
 t^{n+\sigma} = \sigma t^{n+1} + (1-\sigma) t^n, 
 \quad u^{n+\sigma} = \sigma u^{n+1} + (1-\sigma) u^n,
\]  
where $\sigma = \const, \ \sigma \in [0,1]$ is the weight parameter.

For the approximate solution of the problem (\ref{3.2})--(\ref{3.6}), we will use a two-level scheme
\begin{equation}\label{4.1}
 \frac{u^{n+1} - u^n}{\tau} + \sum_{j=1}^{l} \beta_j u^{n+\sigma} + A u^{n+\sigma} - \sum_{j=1}^{l} \beta_j \nu_j w_j^{n+\sigma}
 + \sum_{i=1}^{m} \alpha_i A v_i^{n+\sigma} = \varphi^{n+\sigma} ,
\end{equation} 
\begin{equation}\label{4.2}
 \frac{v_i^{n+1} - v_i^n}{\tau} + \mu_i v_i^{n+\sigma} - u^{n+\sigma} = 0 ,
 \quad i = 1,2, \ldots, m, 
\end{equation} 
\begin{equation}\label{4.3}
 \frac{w_j^{n+1} - w_j^n}{\tau} + \nu_j w_j^{n+\sigma} - u^{n+\sigma} = 0 ,
 \quad j = 1,2, \ldots, l,
 \quad n = 0, 1, \ldots,   
\end{equation}
\begin{equation}\label{4.4}
 u^0 = u_0,
 \quad v_i^0 = 0,
 \quad i = 1,2, \ldots, m, 
 \quad w_j^0 = 0,
 \quad j = 1,2, \ldots, l,  
\end{equation} 
where, for example, $\varphi^{n+\sigma} = f(t^{n+\sigma})$.
The stability of this scheme is established under standard constraints on the weight of $\sigma$.

\begin{thm}\label{t-2}
If (\ref{3.7}) and $\sigma \geq 0.5$ are satisfied, the scheme (\ref{4.1})--(\ref{4.4}) is unconditionally stable and for the approximate solution the estimate 
\begin{equation}\label{4.5}
 \|y^{n+1}\|_* \leq \|u_0\| + 2 \sigma \sum_{k=1}^{n} \tau \|\varphi^{k+\sigma}\| , 
 \quad n = 0, 1, \ldots , 
\end{equation} 
where
\[
 \|y^n\|_*^2 = \|u^n\|^2 + \sum_{i=1}^{m} \alpha_i \|v_i^n\|_A^2 + \sum_{j=1}^{l} \beta_j \nu_j \|w_j^n\|^2 ,
\] 
is true.
\end{thm}

\begin{pf} 
If we multiply equation (\ref{4.1}) by $\tau u^{n+\sigma}$, equation (\ref{4.2}) by $\tau \alpha_i A v_i^{n+\sigma}$, and equation (\ref{4.3}) by $\tau \beta_j \nu_j w_j^{n+\sigma}$, we get the inequality
\begin{equation}\label{4.6}
\begin{split}
 ((u^{n+1} & - u^n),u^{n+\sigma}) + \sum_{i=1}^{m} \alpha_i (A (v_i^{n+1} - v_i^n),v_i^{n+\sigma}) \\
 & + \sum_{j=1}^{l} \beta_j \nu_j ((w_j^{n+1} - w_j^n),w_j^{n+\sigma}) \leq \tau (\varphi^{n+\sigma} ,u^{n+\sigma})
\end{split}
\end{equation} 
Given that
\[
 2 u^{n+\sigma} = (2 \sigma - 1 ) (u^{n+1} - u^n) + (u^{n+1} + u^n) ,
\] 
at $\sigma \geq 0.5$, we have  
\[
 ((u^{n+1} - u^n),u^{n+\sigma}) \geq \frac{1}{2} (\|u^{n+1}\|^2 - \|u^n\|^2) .
\] 
From the inequality (\ref{4.6}), it follows
\begin{equation}\label{4.7}
 \|y^{n+1}\|_*^2 - \|y^n\|_*^2 \leq 2 \tau \||\varphi^{n+\sigma}\| \|u^{n+\sigma}\|. 
\end{equation} 
Under our constraints on $\sigma$, we have
\[
 \|u^{n+\sigma}\| \leq \sigma (\|u^{n+1}\| + \|u^n\|) \leq \sigma (\|y^{n+1}\|_* + \|y^n\|_*) .
\] 
Thus from the inequality (\ref{4.7}),  we obtain the estimate 
\[
 \|y^{n+1}\|_* \leq \|y^n\|_* + 2 \sigma \tau \|\varphi^{n+\sigma}\| .
\] 
From this follows the inequality being proved (\ref{4.5}).  
\end{pf} 

The computational realization of the scheme with weights (\ref{4.1})--(\ref{4.4}) can be organized as follows.
From equations (\ref{4.2}), (\ref{4.3}), we have
\begin{equation}\label{4.8}
\begin{split}
 v_i^{n+1} & = \frac{\sigma \tau }{1 + \sigma \mu_i \tau } u^{n+1} + \eta_i^{n} , \\
 \eta_i^{n} & = \frac{1 }{1 + \sigma \mu_i \tau } \Big ( (1-\sigma) \tau u^n + \big (1 - (1-\sigma) \mu_i \tau \big ) v_i^{n} \Big ) ,
 \quad i = 1,2,\ldots, m ,
\end{split}
\end{equation}  
\begin{equation}\label{4.9}
\begin{split}
 w_j^{n+1} & = \frac{\sigma \tau}{1 + \sigma \nu_j \tau } u^{n+1} + \theta_j^{n} , \\
 \theta_j^{n} & = \frac{1 }{1 + \sigma \nu_j \tau } \Big ( (1-\sigma) \tau u^n + \big (1 - (1-\sigma) \nu_j \tau \big ) w_j^{n} \Big ) ,
 \quad j = 1,2,\ldots, l .
\end{split}
\end{equation} 
Substituting expressions (\ref{4.8}) and (\ref{4.9}) into equation (\ref{4.1}) gives
\begin{equation}\label{4.10}
 (b_I I + \sigma \tau b_A A) u^{n+1} = \chi^n,
\end{equation} 
where $I$ is a unit operator, and
\[
 b_I = 1 + \sigma \tau \sum_{j=1}^{l} \frac{\beta_j}{1 + \sigma \nu_j \tau } ,
 \quad b_A = 1 + \sigma \tau \sum_{i=1}^{m} \frac{\alpha_i}{1 + \sigma \mu_i \tau} .
\] 
For the right-hand side of equation (\ref{4.10}), we have
\[
\begin{split}
  \chi^n & = \Big (1 - (1-\sigma) \tau \big ( \sum_{j=1}^{l} \beta_j +A \big ) \Big ) u^n + \tau \varphi^{n+\sigma} \\
 & + \tau \sum_{j=1}^{l} \beta_j \nu_j \big (\sigma \theta_j^{n} + (1-\sigma)w_j^n \big ) 
   - \tau \sum_{i=1}^{m} \alpha_i A \big (\sigma \eta_i^{n} + (1-\sigma) v_i^n \big ) .
\end{split}
\]
The operator $B = B^* > 0$, where $B = b_I I + \sigma \tau b_A A$, and so equation (\ref{4.10}) is uniquely solvable for $u^{n+1}$.
After finding $u^{n+1}$, the auxiliary quantities $v_i, \ i = 1,2,\ldots, m,$ $w_j, \ j = 1,2,\ldots, l,$ are calculated according to (\ref{4.8}), (\ref{4.9}).  
The increase in the computational complexity of the numerical solution of problems with memory is not fundamental in comparison to the problem without memory effects.

\section{Numerical examples}\label{s:5}

In the above calculations, we have limited ourselves to one-dimensional problems in space.
Considering the current state of numerical solution of boundary value problems for elliptic equations, the transition to more complex multidimensional problems of heat conduction with memory is not of fundamental nature. 
The methodological issues of estimating the rate of convergence of the approximate solution to the exact solution of the system of evolution equations are not discussed. 
We present numerical data on sufficiently detailed computational grids in time and space when we can neglect computational errors.

To illustrate the peculiarities of problems with memory effects, we consider the equation (\ref{3.1}) in which
\[
 \tilde{\alpha}(t) = \alpha \exp(- \mu t),
 \quad \tilde{\beta }(t) = \beta \exp(- \nu t) . 
\] 
i.e., in (\ref{2.6}) $m=l = 1$.
The system of equations (\ref{3.3})--(\ref{3.5}) takes the form
\begin{equation}\label{5.1}
 \frac{d v}{d t} + \mu v - u(t) = 0 ,
\end{equation} 
\begin{equation}\label{5.2}
 \frac{d w}{d t} + \nu w - u(t) = 0 ,
\end{equation} 
\begin{equation}\label{5.3}
 \frac{d u}{d t} + \beta u + A u - \beta \nu w + \alpha A v = f(t) .
\end{equation}  
For (\ref{5.1})--(\ref{5.3}), a Cauchy problem is posed when
\begin{equation}\label{5.4}
 u(0) = u_0,
 \quad v(0) = 0,
 \quad w(0) = 0 .
\end{equation} 

It makes sense to distinguish the following separate classes of heat conduction problems, for which we will give the corresponding equation for temperature.
\begin{description}
 \item[1. Classical thermal conductivity.] The standard model of thermal conductivity without regard to memory effects corresponds to the initial condition (\ref{3.2}) and $\alpha = \beta = 0$:
\[
 \frac{d u}{d t} + A u = f(t) .
\]  
 \item[2. Heat flux memory] In this case $\alpha \neq 0, \ \ \beta = 0$. From (\ref{5.1})--(\ref{5.3}) we get
\[
 \frac{d^2 u}{d t^2} + (\mu + A) \frac{d u}{d t} + (\mu + \alpha) A u = \frac{d f}{d t} + \mu f(t) .
\]   
In this telegraph equation, the term at the first derivative of the solution in time describes the damping effects and the term with the solution --- the oscillatory effects.  
The initial conditions follow from (\ref{5.4}) and equation (\ref{5.3}) at $t=0$: 
\[
 u(0) = u_0,
 \quad \frac{d u}{d t} = f(0) - A u_0 . 
\] 
\item[3. Heat capacity memory] The heat capacity memory effects are accounted for by choosing $\alpha = 0, \ \beta \neq 0$. From (\ref{5.1})--(\ref{5.3}), we get
\[
 \frac{d^2 u}{d t^2} + (\beta + \nu + A) \frac{d u}{d t} + \nu A u = \frac{d f}{d t} + \nu f(t) .
\] 
Again we have a telegraph equation for temperature, for which the initial conditions are given as 
\[
 u(0) = u_0,
 \quad \frac{d u}{d t} = u_1,
 \quad u_1 = f(0) - (\beta + A) u_0 . 
\] 
 \item[4. General case] With $\alpha \neq 0, \ \ \beta \neq 0$ we have a third-order evolutionary equation:
\[
\begin{split}
 \frac{d^3 u}{d t^3} & + (\mu + \nu + \beta + A) \frac{d^2 u}{d t^2} \\
  & + \big ((\mu + \nu)(\beta + A) + \nu (\mu - \beta ) + \alpha A \big ) \frac{d u}{d t}  
  + \nu (\mu + \alpha) A u \\
  & = \frac{d^2 f}{d t^2} + (\mu + \nu) \frac{d f}{d t} + \mu \nu f(t) .
\end{split}
\] 
The initial conditions (\ref{5.4}) on the solutions of the system of equations (\ref{5.1})--(\ref{5.3}) are
\[
 u(0) = u_0,
 \quad \frac{d u}{d t} = u_1 ,
 \quad \frac{d^2 u}{d t^2} = \frac{d f}{d t}(0) - (\beta + A) u_1 + (\beta \nu - \alpha A) u_0 . 
\] 
The system of equations (\ref{5.1})--(\ref{5.3}) can be seen as one possible transition from a third-order evolutionary equation to a first-order system of equations. 
\end{description} 

Numerical illustrations of the features of local generalized thermal conductivity models are given, usually (see, for example, \cite{shen2008notable,hu2009study,zhang2013numerical}) tracing the dynamics of the initial temperature profile. We consider a similar test problem whose numerical solution is carried out with a difference space approximation.
Let us assume that in dimensionless variables $x \in [0,1]$ and a uniform grid with step $h$ is introduced. Let us denote by $\omega$ the set of internal nodes of the grid in space. Given (\ref{2.14}), for $y(x) = 0, \ x \notin \omega$, we define the difference operator $A$ by the relation
\[
 A y = - \frac{1}{h^2} \big (y(x+h) - 2 y(x) + y(x-h) \big ), 
 \quad x \in \omega . 
\] 
The initial condition (\ref{3.2}) is taken as
\[
 u_0(x) = \left \{ \begin{array}{rr}
  x, & 0 < x \leq 0.5 , \\\
  0, & 0.5 < x < 1 , \\
\end{array}
\right .
\quad x \in \omega .
\]  
The calculations were performed on a grid over a space with $h=2\cdot 10^{-3}$, with $0 < t \leq T, \ T = 0.1$, and $\tau = 5\cdot 10^{-5}$.
A purely implicit scheme was used ($\sigma = 1$ in the scheme (\ref{4.1})--(\ref{4.4})).

\begin{figure}
\centering
\includegraphics[width=0.75\linewidth]{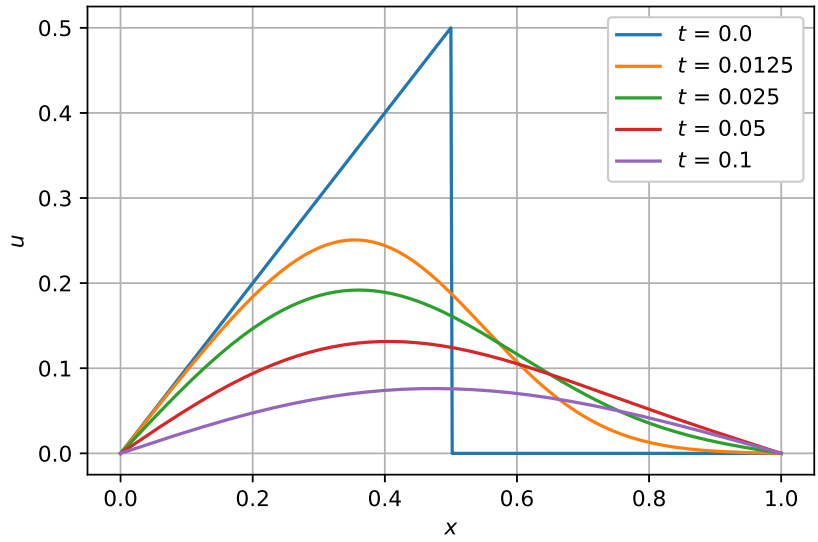} 
\caption{Solution for individual levels in time: model 1 ($\alpha = \beta = 0$).}
\label{f-1}
\end{figure}

\begin{figure}
\centering
\includegraphics[width=0.75\linewidth]{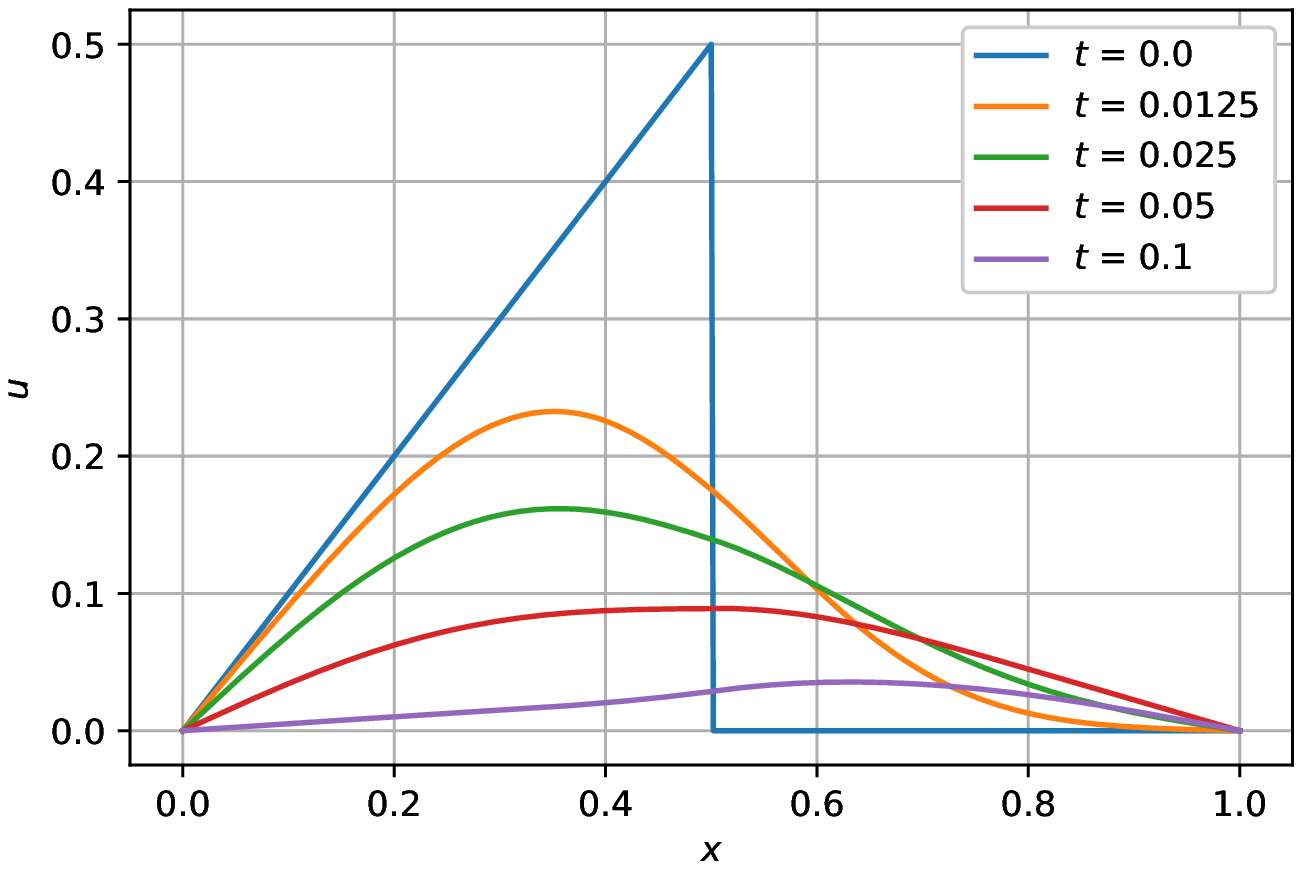} 
\caption{Solution for individual levels in time: model 4 ($\alpha = \beta = 5$, $\mu = \nu  = 1$).}
\label{f-2}
\end{figure}

The dynamics of the process when using model 1 is presented in Fig.~\ref{f-1}, and for model 4 --- in Fig.~\ref{f-2}.  
The effects of heat flux memory and heat capacity lead to a significant rearrangement of the temperature profile.
The influence of individual parameters of the mathematical model is traced to models 2 and 3. 

\begin{figure}
\centering
\includegraphics[width=0.75\linewidth]{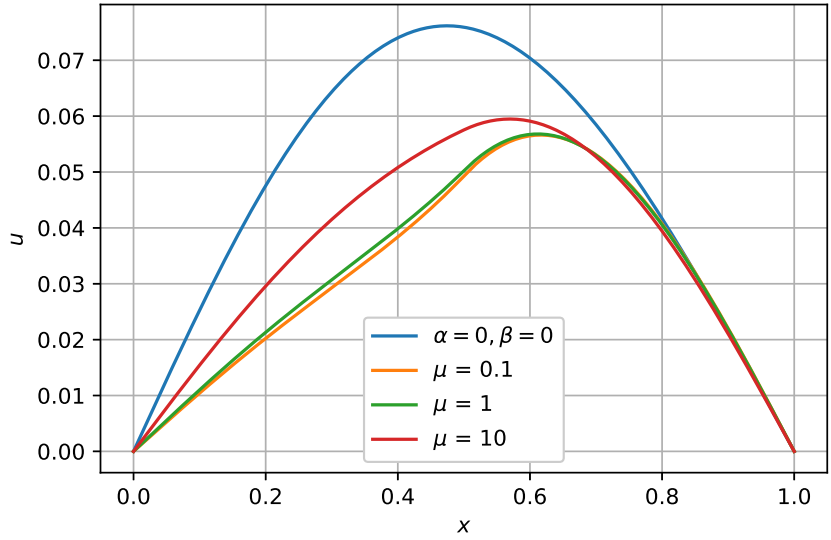} 
\caption{Effect of the parameter $\mu$: solution at $t = 0.1$, model 2 ($\beta = 0$) at $\alpha = 5$.}
\label{f-3}
\end{figure}

\begin{figure}
\centering
\includegraphics[width=0.75\linewidth]{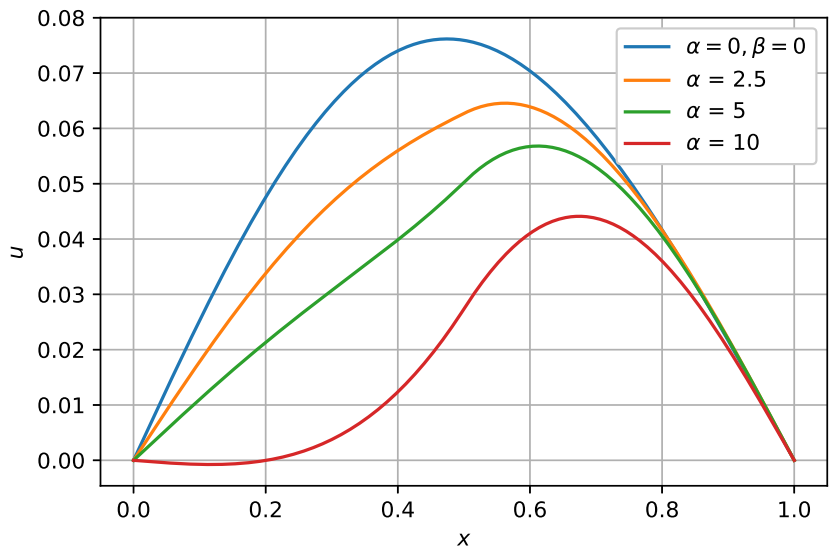} 
\caption{Effect of the parameter $\alpha $: solution at $t = 0.1$, model 2 ($\beta = 0$) at $\mu = 1$.}
\label{f-4}
\end{figure}

The results of heat flow memory (model 2) are presented in Figs.~\ref{f-3}, \ref{f-4}.
The influence of the parameter $\alpha$ is most pronounced. At large $\alpha$ strong oscillatory effects are observed.
In particular, at time $t = 0.1$ at $\alpha=10$, we observe negative values of the solution $u(x,t)$.

\begin{figure}
\centering
\includegraphics[width=0.75\linewidth]{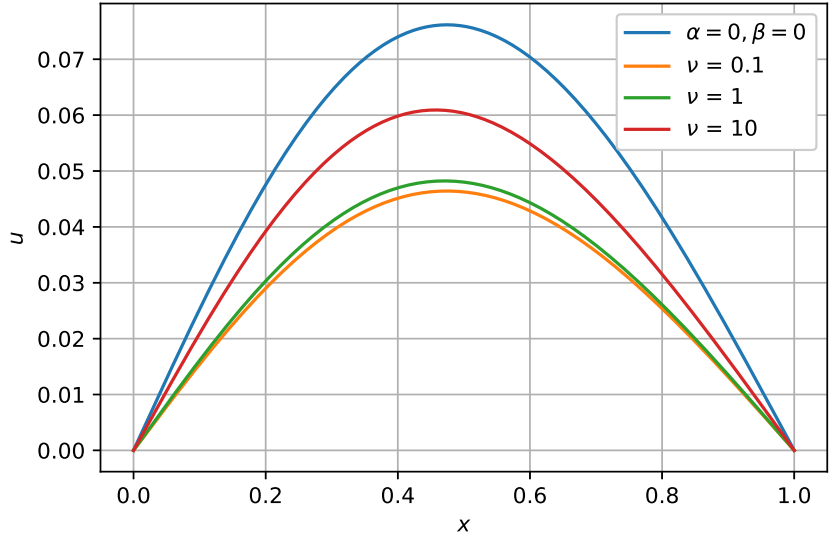} 
\caption{Effect of parameter $\nu$: solution at $t = 0.1$, model 3 ($\alpha = 0$) at $\beta = 5$.}
\label{f-5}
\end{figure}

\begin{figure}
\centering
\includegraphics[width=0.75\linewidth]{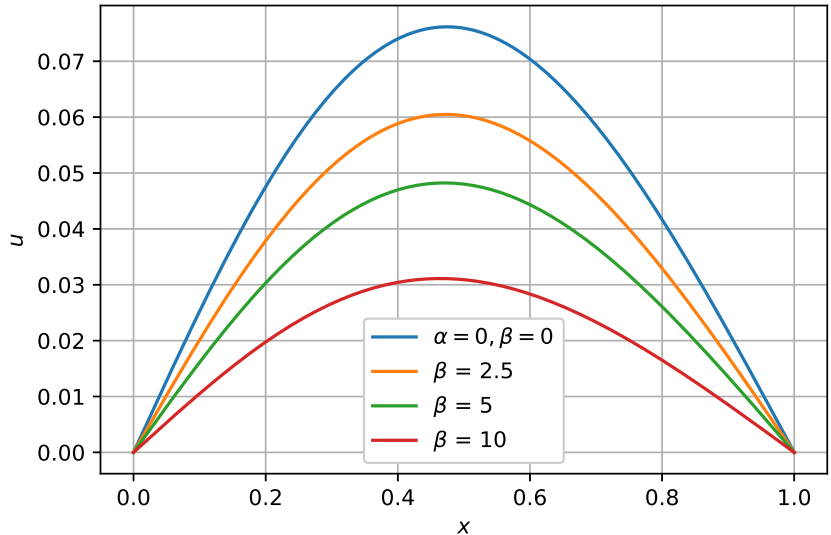} 
\caption{Effect of parameter $\beta$: solution at $t = 0.1$, model 3 ($\alpha = 0$) at $\nu = 1$.}
\label{f-6}
\end{figure}

Similar calculated data for heat capacity memory (model 3) are presented in Figs.~\ref{f-5}, \ref{f-6}.
In this model, the effects of dissipation are most pronounced.

\section{Conclusions}

\begin{enumerate}
 \item A boundary value problem for the first order integrodifferential evolution equation with difference-type kernels, which describes heat conduction processes in solids in linear approximation with allowance for memory effects, is formulated. Processes with a memory of heat flow and memory of heat capacity (internal energy) are modeled.
 \item The nonlocal time problem is reduced to a system of local evolution equations of the first order under the condition that the heat flow and heat capacity relaxation functions are represented as a sum of exponents. We obtained an estimate of the stability of the solution of the Cauchy problem for the system of equations concerning the initial data and the right-hand side in the corresponding Hilbert spaces.
 \item A unconditionally stable scheme with weights for an evolutionary system of equations modeling thermal conductivity in solid media with memory has been proposed and investigated. The transition to a new level in time is not more complicated than the computational implementation of standard schemes with weights for the classical heat conduction equation.
 \item Numerical results of using a purely implicit scheme in approximate solutions of heat conduction problems with allowance for memory effects are presented.
The evolution of the initial state at different values of model parameters is considered.
\end{enumerate} 

\newcommand{\nosort}[1]{}

\end{document}